\documentclass[fleqn]{mat01}
\usepackage{times,mathtimy,amssymb}
\begin{document}

\font\its = cmsl8
\font\rms = cmr7
\font\bfs = cmbx7
\font\gb = cmbx10 scaled \magstep1

\def\sls{\its}
\def\nnn{\hbox{\csc I\kern-.22em N}}
\def\rrr{\hbox{\csc I\kern-.22em R}}
\def\kkk{\hbox{\csc I\kern-.22em K}}
\def\qqq{\hbox{\raise.5ex \hbox{$\scriptscriptstyle |$}\kern-.36em {\csc Q}}}
\def\zzz{\hbox{\bf Z}}

\def\criterion{\trivlist\item[\hskip\labelsep{\it Criterion.}]}
\def\exam{\trivlist\item[\hskip\labelsep{\it Example.}]}
\def\problem{\trivlist\item[\hskip\labelsep{\it Problem.}]}
\newtheorem{theor}{\bf Theorem}
\newtheorem{propo}[theor]{\rm PROPOSITION}
\newtheorem{corol}[theor]{\rm COROLLARY}

\setcounter{page}{225}
\firstpage{225}

\title{Characteristic properties of large subgroups in primary abelian
groups}

\markboth{Peter V~Danchev}{Large subgroups in primary abelian groups}

\author{PETER V~DANCHEV}

\address{Department of Mathematics, Plovdiv University, Paissii
Hilendarski, 4000 Plovdiv, Bulgaria\\
\noindent Current address: 13, General Kutuzov Street, Block 7, Floor 2,
Flat 4, 4003 Plovdiv, Bulgaria\\
\noindent E-mail: pvdanchev@yahoo.com}

\volume{114}

\mon{August}

\parts{3}

\Date{MS received 27 May 2002; revised 19 May 2004}

\begin{abstract}
Suppose $G$ is an arbitrary additively written primary abelian group
with a fixed large subgroup $L$. It is shown that $G$ is (a) summable;
(b) $\sigma$-summable;\break (c)~a $\Sigma$-group; (d) $p^{\omega+1}$-projective
only when so is $L$. These claims extend results of such a kind
obtained by Benabdallah, Eisenstadt, Irwin and Poluianov, {\it Acta
Math. Acad. Sci. Hungaricae} (1970) and Khan, {\it Proc. Indian Acad.
Sci. Sect. A} (1978).
\end{abstract}


\maketitle

\section{Introduction}

The main purpose of this article is to study the relations between the
structures of primary abelian groups and their large subgroups. This is
a semi-expository paper, although some new results have also been
proved. Thus our sections are devoted to producing a number of valuable
examples of sorts of abelian groups where the structure of a $p$-group
$G$ is preserved by its large subgroup $L$ and conversely. The central
theorems here take the form: `$G$ belongs to the class of abelian
$p$-groups ${\cal K}$ if and only if $L$ belongs to ${\cal K}$'.
Evidently this is reminiscent of a number of well-known and documented
statements of the form: `$G$ belongs to ${\cal K}$ if and only if $p^nG$
belongs to ${\cal K}$ for some $n\in {\bf N}$'. Our results, stated in
the present paper, are clearly generalizations of such assertions for
$L$ and in particular for $p^nG$ (see [2,4]). 

In \cite{1}, Benabdallah {\it et~al} proved results of the above type
for the classes of:

\begin{itemize}
\item direct sums of cyclic $p$-groups;
\item direct sums of countable $p$-groups;
\item totally projective $p$-groups;
\item torsion-complete $p$-groups;
\item quasi-closed $p$-groups.
\end{itemize}

Later on, Khan \cite{8} has established such an attainment for
$p$-torsion abelian $\Sigma$-groups but under the additional condition
that they are not separable. We shall supersede here this claim removing
the condition on inseparability, i.e. by proving the attainment for
arbitrary $p$-primary $\Sigma$-groups.

We begin now with an assortment of facts (some of their proofs are quite
elementary, but we include them for the sake of completeness of the
exposition) that motivate this work and that are selected in the next
paragraph; all notations and the terminology used are standard and will
follow essentially the books of Fuchs \cite{6} or those cited in the
bibliography. Nevertheless, for the convenience of readers, we include
certain crucial technical terms and notions. For instance, $p^\omega G =
G^1$ is the intersection of all subgroups $p^n G$ for $n<\omega$. Under
a large subgroup $L$ of $G$ we mean the fully invariant subgroup $L$ of
$G$ with the property $G=L+B$ for every basic subgroup $B$ of $G$. In
particular, $p^n G$ is a large subgroup of $G$ for each natural number
$n$.

We add some additional background material as well from \cite{1} and
\cite{7}, necessary for our successful presentation, because of its
significance to the theory.

If $L$ is an arbitrary large subgroup of a $p$-primary group $G$, then
\cite{10} and [1] independently established the explicit form of a
representation for $L$, that is, $L=\sum_{k<\omega} p^{n_k}G[p^k]$,
where the sequence $n_1,\ldots, n_k,\ldots$ is non-decreasing. From this
it follows by folklore technical arguments that $L^1=G^1$, that
$L[p]=p^mG[p]$ for some $m\in Z^+$ and that $L[p^2]=p^jG[p]+p^tG[p^2]$
for some $j\leq t<\omega$. However, the given representation of $L$ is
not enough to determine whether $L$ inherits the structural
characteristics for $G$ and conversely. That is why another approach
must be demonstrated.

Next, we recall the definition and structural properties (the proofs are
omitted entirely and can be found in \cite{7}) of the so-called
$A$-groups.

If $\mu$ is a limit ordinal, the class $A_{\mu}$ consists of those
$p$-groups $H$ for which there is a containing totally projective
$p$-group $G$ of length not exceeding $\mu$ that satisfies the following
conditions: $H$ is isotype in $G$, \hbox{$p^{\lambda}(G\!/\!H)=
(p^{\lambda}G + H)/\!H$} whenever $\lambda < \mu$ and \hbox{$G\!/\!H$}
is the direct sum of a totally projective group and a divisible group.
The members of the class $A_{\mu}$ are named $\mu$-elementary
$A$-groups. An $A$-group is a direct sum of $\mu$-elementary $A$-groups
for various limit ordinals $\mu$.

Letting $H$ denote an arbitrary reduced abelian $p$-group and $\alpha$
an arbitrary ordinal, $H$ is an $A$-group if and only if both
$p^{\alpha}H$ and \hbox{$H\!/\!p^{\alpha}H$} are $A$-groups.

All of these listed statements will be used extensively with
the corresponding reference.

\section{Invariant properties of large subgroups}

We start with some extensions of totally projective $p$-groups, namely:

\renewcommand\thesubsubsection{\arabic{subsubsection}.}
\subsubsection{\it $p$-torsion $A$-groups.}\hskip -.5pc
Several details on $A$-groups appear in \cite{7}. For example, any
$p$-torsion abelian $A$-group is an isotype subgroup of a totally
projective $p$-group with special properties described in \cite{7}.

\begin{theor}[\!]
Let $G$ be a reduced abelian $p$-group. Then $G$ is an $A$-group if and
only if $L$ is an $A$-group.
\end{theor}

\begin{proof}
First we assume that $G$ is an $A$-group. Therefore \hbox{$G\!/\!G^1$}
is a direct sum of cyclics and $G^1$ is an $A$-group by \cite{7}. But
from \cite{1}, $L^1=G^1$ and hence \hbox{$L\!/\!L^1\leq G\!/\!G^1$.}
Thus $L^1$ is an $A$-group and \hbox{$L\!/\!L^1$} is a direct sum of
cyclics (\cite{6}, p. 110, Theorem~18.1 of L~Kulikov). Finally, again
using \cite{7}, $L$ must be an $A$-group.

Conversely, if $L$ is an $A$-group, then so is $L^1=G^1$ by making use
of \cite{1} and \cite{7}. On the other hand, \hbox{$L\!/\!L^1$} is a direct sum of
cyclics (see \cite{7}) and also it is a large subgroup of
\hbox{$G\!/\!G^1$} owing to \cite{1}. Thus we have in virtue of \cite{1}
that \hbox{$G\!/\!G^1$} is a direct sum of cyclics, and finally \cite{7}
is applicable to get that $G$ is an $A$-group. The proof is
finished.\hfill $\diamondsuit$
\end{proof}

\subsubsection{\it $p$-torsion $C_\lambda$-groups.}\hskip -.5pc
Following Megibben, an abelian $p$-group $C$ is said to be a $C_\lambda$-group if
\hbox{$C\!/\!p^\alpha C$} is totally projective for all ordinals
$\alpha<\lambda$.

Now, we are ready to formulate the following:

\begin{theor}[\!]
Suppose $G$ is an abelian $p$-group and $\lambda$ is an ordinal such
that $\lambda\leq \hbox{ length }(G)$. Then $G$ is a $C_\lambda$-group
if and only if $L$ is a $C_\lambda$-group.
\end{theor}

\begin{proof}
Referring to \cite{1}, we deduce that $p^\alpha G=p^\alpha L$ and
\hbox{$L\!/\!p^\alpha G$} is large in \hbox{$G\!/\!p^\alpha G$} for each
$\alpha\geq \omega$. Consequently, we can apply Theorem~4.5 of \cite{1}
to end the proof. \hfill $\diamondsuit$
\end{proof}

\subsubsection{\it Thin $p$-groups.}\hskip -.5pc
The following definition of a thin $p$-group can be referred to in
\cite{10}. The abelian $p$-group $G$ is thin if every map from a
torsion-complete $p$-group to $G$ is small, that is, the kernel of this
map contains a large subgroup.

We are now in a position to state the following.

\begin{theor}[\!]
Let $G$ be an abelian $p$-group. Then $G$ is thin if and only if $L$ is
thin.
\end{theor}

\begin{proof}
If $G$ is thin, then the same property holds for $L$ as a subgroup (see
\cite{10}).

Assume now that $L$ is thin. Employing \cite{1,6,10}, $G/L$ is a direct
sum of cyclics whence thin. Consequently, taking into account \cite{10},
$G$ would be thin, as desired. The proof is finished. \hfill
$\diamondsuit$
\end{proof}

\subsubsection{\it $p$-torsion $\Sigma$-groups.}\hskip -.5pc
We say that a subgroup of a $p$-torsion abelian group $G$ is high, if it is
maximal with respect to $\cap G^1=0$. There are many non-isomorphic high
subgroups of a given abelian group.

The following definition of a $\Sigma$-group is formulated in \cite{2},
namely, an abelian $p$-group $G$ is termed a $\Sigma$-group provided
that some high subgroup is a direct sum of cyclics. We note that all
high subgroups of a $\Sigma$-group are isomorphic because they are basic
subgroups.

First of all, we list below one useful necessary and sufficient
condition, obtained by us in \cite{3}, for an arbitrary abelian
$p$-primary group to be a $\Sigma$-group.

\begin{criterion}
Suppose $G$ is an abelian $p$-group. Then $G$ is a $\Sigma$-group if and
only if $\displaystyle G[p] = \cup_{n<\omega} G_n\displaystyle$,
$G_n\subseteq G_{n+1}$ and $G_n\cap p^n G = p^\omega G[p]$ for each
natural number $n$.
\end{criterion}
\vspace{.5pc}

We come now to a significant generalisation of \cite{8}.

\begin{theor}[\!]
Suppose $G$ is an abelian $p$-group. Then $G$ is a $\Sigma$-group if and
only if $L$ is a $\Sigma$-group.
\end{theor}

\begin{proof}
Foremost, we recall that $p^\omega L = p^\omega G$, that there is a
natural number $m$ with the property $L[p]=p^m G[p]$ and that $p^n L[p]
= p^{r_n} G[p]$ for each $n<\omega$ and some $m\leq r_n<\omega$ (see
\cite{1}).

Assume now that $G$ is a $\Sigma$-group, i.e., by virtue of the above
criterion, $\displaystyle G[p] = \cup_{n<\omega} G_n\displaystyle$,
$G_n\subseteq G_{n+1}$ and $G_n\cap p^n G=p^\omega G[p]$ for every $n\in
{\bf N}$. Hence $\displaystyle L[p]=\cup_{n<\omega} (G_n\cap
L)\displaystyle$ and $G_n\cap L\subseteq G_{n+1}\cap L$. Moreover, we
calculate $G_n\cap L\cap p^n L = G_n\cap p^n L \subseteq G_n \cap p^n G
= p^\omega G[p] = p^\omega L[p]$. Our criterion leads us to the
conclusion that $L$ is a $\Sigma$-group.

For the reverse inclusion, we assume that $\displaystyle L[p] =
\cup_{n<\omega} L_n\displaystyle$, $L_n \subseteq L_{n+1}$ and $L_n\cap
p^n L = p^\omega L[p]$. Therefore $\displaystyle p^m
G[p]=\cup_{n<\omega} L_n\displaystyle$. Besides, we compute, $L_n\cap
p^{r_n} G[p]= L_n\cap p^n L[p]= p^\omega L[p] = p^\omega G[p]$. Finally,
exploiting our criterion, we detect that $p^m G$ is a $\Sigma$-group,
i.e. the same is valid for $G$, as well (see \cite{2}). The proof is
complete after\break all. \hfill $\diamondsuit$
\end{proof}

\begin{exam}\hskip -.3pc{\it (Megibben).}\ \
The following construction is of a $\Sigma$-$p$-group which is not a
$C_{\omega+1}$-group ($=$ a pillared group in the terminology of Hill).
Let $\displaystyle B=\oplus_{n<\omega} \langle p^n\rangle \displaystyle$
be the direct sum of cyclic groups of order $p^n$ and let $B^{-}$ be its
torsion completion. Then $B$ is pure in $B^{-}$ and \hbox{$B^{-}\!/\!B$}
is isomorphic to $2^{\aleph_0}$ copies of the quasi-cyclic group
$Z(p^{\infty})$. Next, set $\displaystyle
A=\oplus_{\lambda<2^{\aleph_0}} \langle a_{\lambda}\rangle
\displaystyle$, where $\langle a_{\lambda}\rangle \cong \langle p\rangle
$, a cyclic group of order $p$, for each ordinal $\lambda$. If $D$ is
minimal divisible containing $A$, we have \hbox{$D\!/\!A \cong
B^{-}\!/\!B$.} We construct a group $G$ as a subdirect sum of $B^{-}$
and $D$ with kernels $B$ and $A$ respectively. So, with this in hand, it
is not hard to see that $G$ is a $\Sigma$-group with a high subgroup $B$
such that \hbox{$G\!/\!G^1 \cong B^{-}$.}
\end{exam}

\subsubsection{\it Summable $p$-groups.}\hskip -.5pc
The following definition of summable $p$-primary groups is stated in
\cite{6}. We shall say that an abelian $p$-group $G$ of length $\lambda$
($\lambda$ is an ordinal) is summable if $\displaystyle G[p]=
\oplus_{\alpha<\lambda} G_\alpha\displaystyle$ and
\hbox{$G_\alpha\!\setminus\!\{0\}\subseteq p^\alpha G\!\setminus\!
p^{\alpha+1} G$} for all $\alpha<\lambda=\hbox{length}\ (G)$.

\begin{theor}[\!]
Suppose $G$ is a $p$-group. Then $G$ is summable if and only if $L$ is
summable.
\end{theor}

\begin{proof}
First, assume $G$ is summable. Thus $\displaystyle G[p] =
\oplus_{\alpha<\lambda} G_\alpha\displaystyle$, where
\hbox{$G_\alpha\!\setminus\!\{0\}\subseteq p^\alpha$}
\hbox{$G\!\setminus\!p^{\alpha+1} G$} for each $\alpha< \lambda$.
Furthermore, since $L$ is a fully invariant subgroup of $G$ and since
$p^\tau G= p^\tau L$ for all ordinals $\tau\geq\omega$, we conclude
$\displaystyle L[p] = \oplus_{\alpha<\lambda}(G_\alpha\cap
L)\displaystyle$, where \hbox{$(G_\alpha\cap
L)\!\setminus\!\{0\}\subseteq p^\alpha L\!\setminus\!p^{\alpha+1}L$} for
every $\omega\leq\alpha<\lambda$. Moreover, because $p^n L[p] = p^{r_n}
G[p]$ whenever $1\leq n<\omega$; $n\leq r_n<\omega$, we derive that
\hbox{$(G_\alpha\cap L)\!\setminus\!\{0\}\subseteq L\!\setminus\!pL$}
for each $0\leq\alpha<r_1$. Inductively, \hbox{$(G_\alpha\cap
L)\!\setminus\!\{0\}\subseteq pL\!\setminus\!p^2L$} for
$r_1\leq\alpha<r_2,\ldots,$ \hbox{$(G_\alpha\cap
L)\!\setminus\!\{0\}\subseteq p^nL\!\setminus\!p^{n+1}L$} for
$r_n\leq\alpha<r_{n+1}$. So, by setting $\displaystyle
\oplus_{0\leq\alpha<r_1}(G_\alpha\cap L) = L_0\displaystyle$,\ldots,
$\displaystyle \oplus_{r_n\leq\alpha<r_{n+1}}(G_\alpha\cap L) =
L_n\displaystyle$ for every $1\leq n<\omega$ and $G_\alpha\cap
L=L_\alpha$ for every $\alpha\geq\omega$, we infer that $\displaystyle
L[p] = \oplus_{\alpha<\lambda}L_\alpha\displaystyle$, where, for all
ordinals $\alpha<\lambda$, \hbox{$L_\alpha\!\setminus\!\{0\}\subseteq
p^\alpha L\!\setminus\!p^{\alpha+1}L$.} That is why, following the
definition, $L$ is summable.

Now, we assume that $L$ is summable. So, $L^1=G^1$ is summable as its
fully invariant subgroup. On the other hand, conforming with p.~125,
Proposition~84.4 of \cite{6}, $L$ summable implies $L$ is a
$\Sigma$-group. Thus Theorem 4 yields that $G$ is a $\Sigma$-group. For
$H_G$ a high subgroup of $G$ we establish $\displaystyle H_G[p]\oplus
G^1[p]=G[p]\displaystyle$. Since $H_G$ is a direct sum of cyclics, we
may write $\displaystyle H_G[p] = \oplus_{n<\omega} H_n\displaystyle$,
where \hbox{$H_n\!\setminus\! \{0\}\subseteq p^n G\!\setminus\!
p^{n+1}G$} because of the fact that $H_G$ is pure in $G$ (see, for
instance, \cite{6}). Moreover $G^1$ being summable ensures that
$\displaystyle G^1[p]= \oplus_{\alpha<\lambda}C_\alpha\displaystyle$
with \hbox{$C_\alpha\!\setminus\!\{0\}\subseteq p^\alpha
G^1\!\setminus\!p^{\alpha+1} G^1 = p^{\omega+\alpha}
G\!\setminus\!p^{\omega+\alpha+1} G$.} Therefore, we obviously observe
that $\displaystyle G[p]= \oplus_{\beta<\omega} E_\beta \oplus
\oplus_{\omega\leq \beta<\omega+\lambda} E_\beta =
\oplus_{\beta<\omega+\lambda} E_\beta\displaystyle$ by putting $E_n =
H_n$ $\forall n<\omega$ and $E_{\omega+\alpha} = C_\alpha$ $\forall
\alpha<\lambda$. Besides, we calculate that
\hbox{$E_\beta\!\setminus\!\{0\} \subseteq p^\beta
G\!\setminus\!p^{\beta+1} G$} $\forall \beta<\omega+\lambda$. Finally,
we extract that $G$ is summable, as stated. The proof is complete.
\phantom{0000000} \hfill $\diamondsuit$
\end{proof}

\begin{exam}\hskip -.3pc{\it (Hill {\rm \cite{7}}).}\ \
The following is a construction of a summable $C_{\Omega}$-group that is
not a direct sum of countables. For each countable ordinal $\alpha$ let
$G_{\alpha}$ be the unique countable, $p$-primary group of length
$\alpha$ having each of its non-zero Ulm--Kaplansky invariants equal to
$\aleph_0$. Suppose that $\gamma < \Omega$ and that we have embedded
$G_{\alpha}$ in $G_{\beta}$, for each $\alpha<\beta<\gamma$, in such a
way that $p^{\lambda}G_{\beta}\cap G_{\alpha}= p^{\lambda}G_{\alpha}$
holds for all $\lambda$. If $\gamma$ is a limit ordinal, then
$G_{\gamma}\cong \cup_{\alpha<\gamma} G_{\alpha}$ and we have an
embedding of $G_{\alpha}$ in $G_{\gamma}$ such that
$p^{\lambda}G_{\beta}\cap G_{\alpha}= p^{\lambda}G_{\alpha}$ holds for
$\alpha < \beta \leq \gamma$. If $\gamma$ is a non-limit ordinal, we
distinguish two cases.
\end{exam}

\begin{enumerate}
\renewcommand{\labelenumi}{Case \arabic{enumi}.}
\leftskip .6cm
\item $\gamma-1$ is a limit. Let $G_{\gamma-1}$ be a
$p^{\gamma-1}$-high subgroup of $G_{\gamma}$.\vspace{.2pc}

\item $\gamma-2$ exists. Let $G_{\gamma-1}$ be a direct summand
of $G_{\gamma}$.
\end{enumerate}

We remark, in connection with Case 1, that any $p^{\gamma-1}$-high
subgroup of $G_{\gamma}$ is isomorphic to $G_{\gamma-1}$. As far as Case
2 is concerned, $G_{\gamma} \cong G_{\gamma-1}+G_{\gamma}$. Certainly,
it is well-known that $p^{\lambda}G_{\gamma}\cap
G_{\gamma-1}=p^{\lambda}G_{\gamma-1}$ for either type of embedding
described above, so the condition $p^{\lambda}G_{\beta}\cap G_{\alpha}=
p^{\lambda}G_{\alpha}$ continues to hold for $\alpha<\beta\leq \gamma$.

The group of interest is the union ($=$ direct limit) of the groups
$G_{\alpha}$, as $\alpha$ ranges over the countable ordinals, embedded
in one another in the precise manner stated above. Let $G$ denote this
group, namely $G=\cup_{\alpha<\Omega} G_{\alpha}$.

\subsubsection{\it $\sigma$-summable $p$-groups and their direct sums.}
Following Linton and Megibben an abelian $p$-group $C$ is called
$\sigma$-summable if $\displaystyle
C[p]=\cup_{n<\omega}C_n\displaystyle$, $C_n\subseteq C_{n+1}$ and, for
every positive integer $n$, there is an ordinal $\alpha_n$ with $C_n\cap
p^{\alpha_n} C=0$ and $\alpha_n < \hbox{ length } (C)$. Certainly, the
lengths of $\sigma$-summable $p$-groups are infinite ordinals cofinal
with $\omega$.

Next, we concentrate on the following.

\begin{theor}[\!]
Let $G$ be $p$-primary. Then $G$ is $\sigma$-summable if and only if $L$
is $\sigma$-summable.
\end{theor}

\begin{proof}
First, since $G$ is unbounded, $\hbox{ length }(G) = \hbox{ length
}(L)\geq \omega$ (see \cite{1}). Thus, assuming that $G$ is
$\sigma$-summable, we derive that so is $L$ as its subgroup with equal
length, which fact follows immediately from the definition. We only
mention that the case of bounded $G$, hence bounded $L$, follows
directly from \cite{1}.

Conversely, given that $L$ is $\sigma$-summable, hence $\displaystyle
L[p]=\cup_{n<\omega}L_n\displaystyle$, $L_n\subseteq L_{n+1}$ and
$L_n\cap p^{\alpha_n} L=0$ for all $n\geq0$ and some $\alpha_n< \hbox{
length } (G)$. Therefore, as above, $\displaystyle p^m G[p]=
\cup_{n<\omega} L_n\displaystyle$. Since by \cite{1} it holds that
$p^\alpha G= p^\alpha L$ for each ordinal $\alpha\geq\omega$ or that
$p^{\alpha_n} L[p] = p^{k_n} G[p]$ for $\alpha_n<\omega$ and some
natural $k_n\geq \max(\alpha_n, m)$ because in this case $p^{\alpha_n}
L$ is large both in $G$ and $p^{\alpha_n} G$, we obtain $L_n\cap p^{s_n}
G = L_n\cap p^{\alpha_n} L =0$ whenever $s_n< \hbox{ length }(G)= \hbox{
length } (p^m G)\geq \omega$; $s_n=\alpha_n\geq \omega$ or
$\omega>s_n=k_n\geq \max(\alpha_n, m)$.

We construct now subgroups $G_n$ of $G$ in the following manner: $G_n=
\langle a_n \in G[p], L_n | a_n\not\in p^m G \hbox{ and if }$ some
finite sum of degrees of such generating elements $a_n$'s lies in
$p^mG$, then this sum lies in $L_n\rangle$; for more details see \cite{5} too.
Bearing in mind this construction, it is apparent that $\displaystyle
G[p]=\cup_{n<\omega} G_n\displaystyle$, $G_n\subseteq G_{n+1}$.
Moreover, $G_n\cap p^{s_n}G=0$. Indeed, every element of $G_n$ is of the
kind $\varepsilon_1 a_{n1} + \cdots + \varepsilon_s a_{ns} + b_n$, where
$\varepsilon_1, \ldots,\varepsilon_n \geq 0$ are integers; $a_{n1},
\ldots, a_{ns}$ are generators and $b_n\in L_n$. If $\varepsilon_1
a_{n1} + \cdots + \varepsilon_s a_{ns}\in p^m G$, then by definition it
is in $L_n$, so by the above computations the studied intersection is
precisely 0. In the remaining case when $\varepsilon_1 a_{n1} + \cdots +
\varepsilon_s a_{ns}\not\in p^m G$, we find that $\varepsilon_1 a_{n1} +
\cdots + \varepsilon_s a_{ns} +b_n\not\in p^m G$ since $b_n\in
L_n\subseteq p^m G$, i.e. it has $\hbox{ height } <m\leq s_n$, hence the
considered intersection is again equal to 0.

This gives the desired claim and the proof is finished. \hfill
$\diamondsuit$
\end{proof}

\begin{propo}$\left.\right.$\vspace{.5pc}

\noindent Let $G$ be a direct sum of $\sigma$-summable abelian
$p$-groups. Then so is $L$.
\end{propo}

\begin{proof}
Write down $\displaystyle G=\oplus_{i\in I} G_i\displaystyle$, where all
summands $G_i$ are $\sigma$-summable groups. Therefore $\displaystyle
L=\oplus_{i\in I} (L\cap G_i)\displaystyle$, because it is fully
invariant in $G$. But since all components $G_i$ are isotype in $G$,
invoking \cite{1}, we deduce that $L\cap G_i$ are large in $G_i$. By
what we have just argued in Theorem 6, $L\cap G_i$ are $\sigma$-summable
groups, too. This verifies the statement, thus finishing the proof.
\hfill $\diamondsuit$
\end{proof}

\subsubsection{\it $p^{\omega+1}$-projective $p$-groups.}\hskip -.5pc
First of all, we recall the criterion for $p^{\omega+n}$-projectivity
obtained by Nunke \cite{9} when $n\geq 0$ is an integer.

\begin{criterion}\hskip -.4pc \cite{9}.\ \ The abelian $p$-group $G$ is
$p^{\omega+n}$-projective if and only if there exists a subgroup $E \leq
G[p^n]$ such that \hbox{$G\!/\!E$} is a direct sum of cyclics.
\end{criterion}

Clearly, subgroups of $p^{\omega+n}$-projective abelian $p$-groups are
themselves $p^{\omega+n}$-projective.

The following slight modification of a classical result due to
Dieudonn\'e was proved by us in \cite{5}; actually we have established a
more general version for a larger class of abelian groups. Here we
emphasize it, and freely use it in the sequel, without further\break comments.

\begin{criterion}\hskip -.4pc (Dieudonn\'e, 1952).\ \
Suppose $G$ is an abelian $p$-group for which there exists a subgroup
$A$ such that \hbox{$G\!/\!A$} is a direct sum of cyclics. Then $G$ is a
direct sum of cyclics if and only if $\displaystyle A[p]=\cup_{n<\omega}
A_n$, $A_n \subseteq A_{n+1}$ and $A_n\cap {p^n}G = 0$ for every
non-negative\break integer $n$.
\end{criterion}

Major consequences that parallel the example listed below are the
following:

\begin{enumerate}
\renewcommand{\labelenumi}{\arabic{enumi})}
\item If $G$ is a separable $p$-group such that there is $A \leq G$ with
$|A|=\aleph_0$ and \hbox{$G\!/\!A$} is a direct sum of cyclics, then $G$ is
a direct sum of cyclics.

\hskip 1pc This holds easily by making use of the Dieudonn\'e criterion
since $A=\cup_{n<\omega} A_n$, $A_n\leq A_{n+1}$ and all $A_n$ are
finite whence bounded in $G$.

\item (A result due to Fuchs, Mostowski--Sosiada; \cite{6}, p.~111,
Proposition 18.3.) If $G$ is a $p$-torsion group so that there exists $A
\leq G$ with $A$ a direct sum of cyclics and \hbox{$G\!/\!A$} bounded,
then $G$ is a direct sum of cyclics.

\hskip 1pc This follows by the usage of the Dieudonn\'e criterion
because of the fact that there is a natural number $t$ with the property
$p^tG\subseteq A=\cup_{n<\omega} A_n$, $A_n\subseteq A_{n+1}$ and
$A_n\cap p^{t+n}G\subseteq A_n\cap p^nA=0$.
\end{enumerate}

We now begin with the following:

\begin{theor}[\!]
Suppose $G$ is an abelian $p$-group. Then $G$ is
$p^{\omega+1}$-projective if and only if $L$ is
$p^{\omega+1}$-projective.
\end{theor}

\begin{proof}
As we have observed above, $G$ being $p^{\omega+1}$-projective yields
that so does $L$.

To treat the converse, in view of Nunke's ctiterion stated above, we may
write down \hbox{$(L\!/\!E)[p] = \cup_{n<\omega} (L_n\!/\!E)$,}
$L_n\subseteq L_{n+1} \leq L$ and, for each $n<\omega$, $L_n\cap p^n
L\subseteq E$ for some $E\leq L[p]$. Certainly $L_n \subseteq L[p^2]$,
$\forall n<\omega$. Since $p^n L[p^2] = p^{t_n} G[p^2]+ p^{j_n} G[p]$
for some $n \le j_n \le t_n<\omega$ (see \cite{1}), hence $p^{t_n}
G[p^2] \subseteq p^n L[p^2] \subseteq p^{j_n} G[p^2]$, we have $L_n \cap
p^{t_n} G = L_n \cap p^{t_n} G[p^2] \subseteq L_n \cap p^n L[p^2] = L_n
\cap p^nL \subseteq E$. Thus the heights of elements from
\hbox{$L_n\!/\!E$} are bounded in \hbox{$G\!/\!E$} for all $n<\omega$.

Moreover, according to \cite{1}, \hbox{$(G\!/\!E)/(L\!/\!E)\cong
G\!/\!L$} is a direct sum of cyclics. Finally, in the spirit of the
modified variant of the classical Dieudonn\'e criterion formulated
above, that is an expansion of the classical Kulikov's criterion
\cite{6}, we extract that \hbox{$G\!/\!E$} is a direct sum of cyclics.
That is why, in conjunction with the Nunke's criterion, $G$ is
$p^{\omega+1}$-projective, as asserted. The proof is over.
\hfill$\diamondsuit$
\end{proof}

\begin{exam}\hskip -.4pc {\it (Dieudonn\'e).}\ \ The following ingenious
example, due to Dieudonn\'e, demonstrates that there is a separable
$p^{\omega+1}$-projective group that is not a direct sum of cyclics.
Given that $C=\prod_{1\leq k<\omega} \langle c_k \rangle$, where order
$(c_k)=p$ and $C_n=\prod_{n<k<\omega} \langle c_k \rangle$ $\forall
n\geq 1$. Clearly $pC=0$, hence $C$ is a bounded direct sum of cyclic
groups. For every $0\not= x \in c_n +C_n$, let us define the elements
$a_{ni}=a_{ni}(x)$ such that $p^na_{ni}=x$ whenever $i\in I$ for an
index set $I$. Thereby, we construct $K=\langle a_{ni}, C|n\geq 1, i\in
I\rangle = \langle a_{ni}|n\geq 1, i\in I \rangle + C$. Evidently,
$a_{ni}\notin C$ since $pa_{ni}\not= 0$, and moreover
$|C|=2^{\aleph_0}>\aleph_0$ and $|I|\geq |C|$. Therefore $|K|=|I|$.
\end{exam}

First of all, we quote with some comments the following properties of
the so-constructed group $K$ (see also \cite{6}, p.~16, Exercise 11):

\begin{enumerate}
\renewcommand{\labelenumi}{(\alph{enumi})}
\item order $(a_{ni})=p^{n+1}$. For fixed $n$ and $x$, the equation
$p^ny=x$ has at least $2^{\aleph_0}$ solutions, so the index $i$ runs on
$I$ with $|I|\geq 2^{\aleph_0}$. For instance, $a_{nj}=a_{ni}+a_{n-1,i}$
for $i\not=j$ is a guarantor for this, namely $p^na_{nj}=p^na_{ni}$ and
order $(a_{nj})= \hbox{order}\ (a_{ni})$. Henceforth, $\langle
a_{ni}\rangle \cap \langle a_{nj} \rangle = \langle p^na_{ni} \rangle =
p^n \langle a_{ni} \rangle$. Moreover, it is self-evident that $\langle
a_{ni}\rangle \cap C = \langle p^na_{ni} \rangle$. In that aspect,
$\langle a_{ni} \rangle \cap \langle a_{n+1,j} \rangle =0$ $\forall$
$i\not= j$, and $\langle a_{ni}(x)\rangle \cap \langle a_{ni}(y)\rangle
=0$ $\forall$ $x\not= y\in c_n+C_n$ because $p^na_{ni}(x) \not=
p^na_{nj}(y)$.

\item \ \!\hbox{$K\!/\!C$} is an unbounded and uncountable direct sum of
cyclics. Really, $K\!/\!C=\langle a_{ni}|n\geq 1, i\in I\rangle +C\!/\!C
\cong \langle a_{ni}|n\geq 1, i\in I\rangle/\langle a_{ni}|n\geq 1, i\in
I\rangle \cap C = \langle a_{ni}|n\geq 1, i\in I\rangle/p^n\langle
a_{ni}|n\geq 1, i\in I\rangle \cong \oplus_{n\geq 1} \oplus_{i\in I}
[\langle a_{ni}\rangle/p^n\langle a_{ni}\rangle]$, since $[(\langle
a_{ni}\rangle +C)/C]\cap [(\langle a_{nj}\rangle +C)/C]=[(\langle
a_{ni}\rangle +C)\cap (\langle a_{nj}\rangle +C)]/C = [(\langle
a_{ni}\rangle \cap \langle a_{nj}\rangle)+C]/C=(\langle p^na_{ni}\rangle
+C)/C=0$.

\item For each $0\not= x \in c_n+C_n$ it follows that \hbox{$x\in
p^nK\!\setminus\!p^{n+1}K, x\in C\!\setminus\!pC$} and $p^{\omega}K=0$.

\item For every $S\leq K$ such that $S$ is with elements of heights
bounded in $K$ it holds that $|S\cap C|< \aleph_0$.
\end{enumerate}

With the aid of the above Nunke's criterion and (d) it follows at once
that $K$ is $p^{\omega+1}$-projective but not $p^{\omega}$-projective
since otherwise $K=\cup_{n<\omega} K_n$, $K_n\subseteq K_{n+1}$ and
$K_n\cap p^nK=0$ yield $C=\cup_{n<\omega} (K_n\cap C)$ is countable,
which is obviously wrong. We indicate that $C\not= K[p]$ and that
$C\not= pK$.

Next, we shall find some additional properties of basic and large
subgroups in $K$. In fact, $pK=\langle pa_{ni}|n\geq 1, i\in I\rangle$
is a large subgroup which is not a direct sum of cyclics. For a basic
subgroup $B$ we have $B=\oplus_{n\geq 1} \oplus_{i\in I} \langle b_{ni}
\rangle$, where order$(b_{ni})=p^{n+1}$ and $b_{ni}$ are generators
which depend on $a_{ni}$; however we omit the technical details of their
explicit representation. By definition, $B$ is pure in $K$, and
\hbox{$K\!/\!B$} is divisible that is $K=pK+B$. We also observe that $B$ is
unbounded since $K$ is separable unbounded. Moreover, $K$ being
$p^{\omega+1}$-projective implies that it is fully starred, hence
starred. Thus, $|K|=|B|\geq 2^{\aleph_0}$. Besides, $|B\cap C|=\aleph_0$.
Indeed, writing $B_n=\oplus_{1\leq k \leq n} \oplus_{i \in I} \langle
b_{ki}\rangle $, we establish $B_n\subseteq B_{n+1}$, $B_n\cap
p^nB=B_n\cap p^nK =0$ and $B=\cup_{n<\omega} B_n$. In lieu of (d), for
each $n\geq 1$ we obtain $B_n\cap C$ is finite, so we are done.

\begin{problem}
Suppose $C\leq A$ such that both $C$ and \hbox{$A\!/\!C$} are direct
sums of $p$-cyclic groups. Then what is the algebraic structure of $A$?
\end{problem}
\vspace{.5pc}

We know that $A$ is fully starred (by Irwin--Richman), but however a
more precise description is needed. For example, if $A$ is a $p$-group
and $C$ is bounded at $p^n$, we have known via the Nunke's criterion
that $A$ is $p^{\omega+n}$-projective.

\subsubsection{\it Pure-complete $p$-groups.}\hskip -.5pc
The following definition of a pure-complete $p$-group is stated in
\cite{6}. An abelian $p$-group $G$ is pure-complete if for every $S\leq
G[p]$ there is a pure subgroup $T$ of $G$ so that $S=T[p]$.

\begin{propo}$\left.\right.$\vspace{.5pc}

\noindent Suppose $G$ is an abelian $p$-group. If $G$ is
pure-complete{\rm ,} then so is $L$.
\end{propo}

\begin{proof}
Assume that $S\leq L[p]$. Hence $S\leq G[p]$ and $S=M[p]$ for some pure
subgroup $M$ of $G$. An application of \cite{1} assures that $M\cap L$
is large in $M$, whence it is only a routine matter to check that $M\cap
L$ is pure in $L$. Finally, we may write $S=M[p]\cap L[p]=(M\cap L)[p]$,
i.e. by definition $L$ is pure-complete, as claimed. The proof is
completed. \hfill $\diamondsuit$
\end{proof}

\begin{corol}$\left.\right.$\vspace{.5pc}

\noindent Let $G$ be an abelian $p$-group. Then $G$ is pure-complete if
and only if $p^nG$ is pure-complete for some arbitrary fixed positive
integer $n$.
\end{corol}

\begin{proof}
The necessity follows automatically from the previous proposition via
the substitution $p^nG=L$.

After this, the more difficult reverse question is the focus of our
examination. In fact, since $p^nG=p(p^{n-1}G)$, a routine transfinite
induction reveals that we may presume $n=1$. Now, given a subgroup
$S\subseteq G[p]$ whence $S\cap pG\subseteq (pG)[p]$ and $S\cap pG=P[p]$
for some pure subgroup $P$ of $pG$. By referring \cite{2}, we infer that
there is a pure subgroup $K$ of $G$ so that $pK=P$ and
$K[p]=P[p]=(pK)[p]$. Therefore, $S\cap pG=(pK)[p]$. Our future aim,
which we pursue, is to check that $S=T[p]$ for some pure subgroup $T$ of
$G$. In fact, we construct the wanted group $T$ as follows:
\hbox{$T=K+(S\!\setminus\!pG)\cup\{0\}$.} We need to check that $T$ is
indeed a group. In fact, it is self-evident that $0\in T$ and that if
$u\in T$ then $-u\in T$. Next, let $u\in T$ and $v\in T$. Thus $u = a+b$
for $a\in K$ and \hbox{$b\in (S\!\setminus\!pG)\cup\{0\}$} as well as $v
= c+d$ for $c\in K$ and \hbox{$d\in (S\!\setminus\!pG)\cup\{0\}$.}
Therefore $u+v = (a+c) + (b+d)$. It is trivial that $a+c\in K$. Further,
we distinguish two basic cases. First, if \hbox{$b+d\in
(S\!\setminus\!pG)\cup\{0\}$,} there is nothing to prove. Otherwise,
when $b+d\in pG$, we detect that $b+d\in S\cap pG = P[p] = K[p]$ whence
$u+v\in K\subseteq T$. So, in both situations, $u+v\in T$, as needed.
Besides, \hbox{$T[p]=(K+(S\!\setminus\!pG)\cup\{0\})[p]=K[p]+
(S\!\setminus\!pG)\cup\{0\}$} because
\hbox{$(S\!\setminus\!pG)\cup\{0\}= ((S\!\setminus\!pG)\cup\{0\})[p]$.}
Henceforth \hbox{$T[p]= P[p]+ (S\!\setminus\!pG)\cup\{0\}\subseteq S$}
and \hbox{$S=(S\cap pG)\cup (S\!\setminus\!pG)=K[p] \cup
((S\!\setminus\!pG)\cup\{0\})[p]= (K\cup (S\!\setminus\!pG)\cup\{0\})[p]$}
\hbox{$\subseteq (K+ (S\!\setminus\!pG)\cup\{0\})[p]=T[p]$.} Thereby $S=T[p]$.

Finally, we shall verify that $T$ is pure in $G$. Indeed, motivated by
\cite{6}, for each natural number $t$ we compute \hbox{$T[p]\cap p^tG =
S\cap p^tG = [K[p]+(S\!\setminus\!pG)\cup\{0\}]\cap p^tG = $}
\hbox{$[(pK)[p]+(S\!\setminus\!pG)\cup\{0\}]\cap p^tG = K[p]\cap p^tG =
(p^tK)[p] = (p^tT)[p]$.} This concludes the proof. \hfill $\diamondsuit$
\end{proof}

\begin{problem}
If $U$ is a pure subgroup of $L$, then does there exist a pure
subgroup $V$ of $p^nG$ for some non-negative integer $n$ such that $V[p]=U[p]$.
\end{problem}

If the answer is yes, $L$ being pure-complete does imply that the same
is $p^nG$ for this $n$.

\section{Unanswered questions}

Here we quote some major problems that remain unsolved:
\noindent Is it true that $G$ is a (1) direct sum of $\sigma$-summable
$p$-groups;
(2) $p^{\omega+n}$ projective $p$-group;
(3) thick $p$-group;
(4) essentially finitely indecomposable $p$-group;
(5) pure-complete $p$-group;
(6) semi-complete $p$-group;
(7) direct sum of torsion-complete $p$-groups (\cite{1}, Problem 4);
(8) Fuchs 5 $p$-group;
(9) Q-$p$-group;
(10) $\aleph_1$-separable $p$-group;
(11) weakly $\aleph_1$-separable $p$-group
if and only if so is $L$?

\looseness 1 However, not all properties of $G$ are inherited for $L$
and conversely. In fact, we give the following example that is our goal
here: The abelian $p$-group $G$ is said to be starred (see, for
instance, \cite{10}) if it has the same power as its basic subgroup $B$.
Set $L=p^n G$. Suppose now that $G$ is unbounded torsion-complete with a
basic subgroup $\displaystyle B=\oplus_{0<k\leq n} \oplus_{\aleph_1}
Z(p^k) \oplus\oplus_{n<k<\omega} Z(p^k)\displaystyle$ of cardinality
$\aleph_1$, and assume that the {\it generalized continuum hypothesis}
holds. Then we claim that $G$ is starred, and $p^n G$ is unbounded
torsion-complete (see, for example, \cite{2,6}) but not starred. This is
fulfilled because $p^n B$ is a basic subgroup of $p^n G$ (see \cite{6})
and so complying with p.~29, Exercise~7 of \cite{6} we compute $|G| =
|B|^{\aleph_0} = \aleph_1^{\aleph_0} = \aleph_0^{\aleph_0} = \aleph_1 =
|B| = |p^n G| = |p^n B| ^{\aleph_0} > \aleph_0 = |p^n B|$, that
substantiates our claim.

A similar example can be given to show that $|G|>|L|$ in general. Thus
the assumption that $L$ is countable does not imply that $G$ is
countable.

\section*{Acknowledgements}

The author is very grateful to the expert referee for helpful
suggestions necessary for the careful preparation of the present
manuscript.

\end{document}